\documentclass[reqno]{amsart}
\usepackage{amsmath}
\usepackage{amssymb}
\usepackage{stmaryrd}
\allowdisplaybreaks[4]
\title{
The Number of Hypergraphs and Colored Hypergraphs with 
Hereditary Properties
}
\author{Yoshiyas Ishigami}
\address{Department of Information and Communication Engineering, 
The University of Electro-Communications, Chofu, Tokyo 182-8585, Japan.
}
\email{yoshiyas@ice.uec.ac.jp}
\thanks{
This research was partially supported by the Ministry of Education, 
Culture, Sports, Science and Technology 
of Japan, Grant-in-Aid for Scientific Research (C)(19540119).
}
\date{\today}
\def\picture #1 by #2 (#3){
                \vbox to #2{
                        \hrule width #1  height 0pt depth
0pt
                        \vfill
                        \special{picture #3}}}
\def\scaledpicture #1 by #2(#3 scaled #4){{
                \dimen0=#1 \dimen1=#2
                \divide\dimen0 by 1000 \multiply\dimen0 by
#4
                \divide\dimen1 by 1000 \multiply\dimen1 by
#4
                \picture\dimen0 by \dimen1 (#3 scaled #4)}}
\newtheorem{adf}{Definition}[section]
\newtheorem{tha}{Theorem}[section]

\newtheorem{thm}{Theorem}[section]

\newtheorem{claim}[thm]{Claim}

\newtheorem{remark0}[thm]{Remark}

\newtheorem{algorithm0}[thm]{Algorithm}
\newtheorem{asett}[thm]{Setup}

\newenvironment{df}{
\begin{adf}\begin{sl}}
{\end{sl}\thqed\end{adf}}
\newenvironment{sett}{
\begin{asett}\begin{sl}}{\end{sl}\thqed\end{asett}}

\newcommand{\thqed}{\hfill\fbox{}\\ }

\newcommand{\claimqed}{\hfill\bigskip\fbox{}\\}

\def\brkt#1{\left({#1}\right)}

\def\ang#1{\langle{#1}\rangle}

\newcommand{\Prob}{{\mathbb P}}
\newcommand{\Ex}{{\mathbb E}}

\newtheorem{procedure0}{Procedure}

\newcount\timehh\newcount\timemm\timehh=\time\divide\timehh by 60
\timemm=\time\count255=\timehh\multiply\count255 by-60
\advance\timemm by \count255
\ifnum\timehh=12\def\apm{pm}\else
\ifnum\timehh>12\def\apm{pm}\advance\timehh by-12\else
\def\apm{am}\fi\fi
\def\timestamp{\number\timehh\,:\,\ifnum\timemm<10 0\fi\number\timemm\,\apm}
\begin{document}
\maketitle
\begin{abstract}
Using Szemer\'edi's regularity lemma, 
Erd\H{o}s, Frankl and R\"{o}dl (1986) \cite{EFR86} showed that 
for any monotone family of graphs ${\mathcal P}$, the number of 
graphs on vertex set $[n]$ in ${\mathcal P}$ is 
$$ 2^{(1+o(1))
{\rm ex}(n,{\mathcal P}){n\choose 2}
}
$$ where ${\rm ex}(n, {\mathcal P}){n\choose 2}$ is the maximum number of 
edges of an $n$-vertex graph which has no edge-induced subgraph in 
${\mathcal P}.$
It was extended from monotone families to hereditary families, by 
Alekseev (1993) \cite{A93} and by Bollob\'as and Thomason (1997) \cite{BT97}. 
Kohayakawa et al. (2003) \cite{KNR} further extended it from graphs 
to $3$-uniform hypergraphs, using Frankl-R\"{o}dl(2002)'s version of 
$3$-uniform hypergraph regularity lemma.
We will extend it to $k$-uniform hypergraphs and to $k$-uniform 
colored hypergraphs.
Our proof may be a simple example illustrating how to apply a new 
hypergraph regularity lemma by \cite{I06}.
\end{abstract}
\section{Introduction}
\subsection{Notation and statement of the main result}
Given a positive integer $k$ and a set $V,$ 
we denote 
${V\choose k}:=\{e\subset V; |e|=k\}.
$
For a set $C,$ 
a function $H: {V\choose k}\to C$
 is a 
{\bf  $k$-uniform $C$-colored hypergraph} 
(or {\bf $(k,C)$-graph}) 
where 
members in $V$, in ${V\choose k}$, and in $C$ are 
called {\bf vertices}, {\bf edges}, and {\bf colors.}
The sets $V$ and $C$ are  called a {\bf vertex set} and a {\bf color set}.
When $H$ is a $(k,C)$-graph, 
$V(H)$ means the vertex set of $H$ and 
if $|V(H)|=n$, $H$ is called to be {\bf $n$-vertex}.
When $C=\{black, white\},$ 
an ($n$-vertex $k$-uniform $C$-colored) 
hypergraph is considered to be an ordianary 
($n$-vertex $k$-uniform) hypergraph, which have been studied by many researchers.
When the important information on the color set $C$ 
is often $|C|$ only, sometimes a $(k,C)$-graph 
is called simply a $k$-uniform {\bf $|C|$-colored } hypergraph 
or a {\bf $(k,|C|)$-graph}.
A {\bf subgraph} of a $(k,C)$-graph $H$ is a $(k,C)$-graph 
obtained from $H$ by deleting some vertices of $H$ (if necessary).
\par 
A (finite or infinite) family of $k$-uniform $C$-colored hypergraphs 
is said to be a {\bf $(k,C)$-property}, (or a {\bf $(k,|C|)$-property} 
or simply a {\bf property}), when 
if the family contains a $(k,C)$-graph, say $H$, 
then the family also contains any $(k,C)$-graph $H'$
obtained from $H$ by relabeling the vertices of $H.$ 
(That is, the labels (i.e. the names) of the vertices are 
irrelevant for the property. However we distinguish the labels of the colors, 
so a hypergraph with all edges black and a hypergraph with all 
edges white are considered to be different.
)
Note that when two $(k,C)$-graphs $H_1,H_2$ in the same property,
the numbers of the vertices in $H_1$ and $H_2$ are not necessarily equal.
A $(k,C)$-graph {\bf satisfies the property} 
if-and-only-if it belongs to the property.
When ${\mathcal P}$ is a $(k,C)$-property, we usually denote by ${\mathcal P}_n$ 
the family of the $n$-vertex $(k,C)$-graphs 
on vertex set $[n]:=\{1,2,\cdots, n\}$ 
satisfying ${\mathcal P}$.
Clearly  $$|{\mathcal P}_n|\le |C|^{n\choose k}.$$
And a property is {\bf hereditary } if-and-only-if, 
whenever a colored hypergraph satisfies the property, any subgraph 
of it also satisfies 
the property.
\par
When 
${\mathcal P}$ is a $(k,C)$-property, 
we say that an $n$-vertex $(k, 
2^C\setminus \{\emptyset\})$-graph $H$ 
is {\bf ${\mathcal P}$-good} 
if-and-only-if 
 ${\mathcal P}$ contains 
any $n$-vertex $(k,C)$-graph $H'$ obtained from $H$ 
by recoloring each edge $e\in {V(H)\choose k}$ 
with {\em any} member $c\in H(e)\subset C$.
\par
For a $(k,C)$-property, we define 
\begin{eqnarray}
{\rm ex}(n, {\mathcal P}):=\max_H\{ 
\Ex_{e\in {[n]\choose k}}
\log_2 |H(e)|
; H
\mbox{ is a 
${\mathcal P}$-good 
$(k,2^C\setminus\{\emptyset\})$-graph on vertex set $[n]$
 } 
\}
\end{eqnarray}
where $\Ex$ means the expectation or average, 
i.e. $\Ex_e={1\over {n\choose k}}\sum_e$ in the above.
\begin{thm}[Main Theorem]\label{uqq03}
Let $k$ be a positive integer and a finite set $C.$
When $k=O(1)$ and $|C|=O(1)$ as $n\to\infty,$ if 
${\mathcal P}$ is a $(k,C)$-property 
then the number $|{\mathcal P}_n|$ of $(k,C)$-graphs on vertex set 
$[n]=\{1,2,\cdots,n\}$ satisfying ${\mathcal P}$ is 
$$|{\mathcal P}_n|= 2^{({\rm ex}(n,{\mathcal P})+o(1)){n\choose k}}.$$
\end{thm}
\subsection{Basic remarks}
In our main theorem, the $\ge$-part will be easily seen, so 
the $\le$-part is the main part of our result.
Our proof is constructive, so $o(1)$ is bounded by a certain 
function, though we will not write the explicit form of the function.
\par
It is easily seen that 
$\lim_{n\to\infty} {\rm ex}(n, {\mathcal P})$ 
exists, because an easy averaging argument implies that 
${ {\rm ex(n, {\mathcal P})}}$ is 
non-increasing for $n.$ (The argument will be seen in the early part of 
the proof of the theorem.)
Thus by our main theorem, $\lim_{n\to\infty}{\log |{\mathcal P}_n|\over 
{n\choose k}}$ exists.
\par
Given a $(k,C)$-property ${\mathcal F},$ 
we denote by ${\rm Forb}({\mathcal F})$ 
the $(k,C)$-property which consists of all 
$(k,C)$-graphs containing no copy of $F$ 
as a subgraph for any $F\in {\mathcal F}.$ 
Also denote by $
{\rm Forb}(n, {\mathcal F})
={\rm Forb}(n, {\mathcal F})
$ 
the family of such hypergraphs on vertex set $[n].$
It is easy to see that ${\rm Forb}({\mathcal F})$ is 
always hereditary. On the other hand, 
any hereditary $(k,{C})$-property can be expressed in this way.
In fact, for any hereditary $(k,C)$-property ${\mathcal P}$, 
let ${\mathcal F}$ be the $(k,C)$-graphs (with any number of vertices) 
which does not satisfy ${\mathcal P}.$ This family and its members 
are called {\bf forbidden} for 
${\mathcal P}.$
Then it is easily seen that ${\mathcal P}={\rm Forb}({\mathcal F})$ 
and ${\mathcal P}_n={\rm Forb}_{\rm }(n,{\mathcal F})$.
(Indeed, if ${H}$ does not satisfy ${\mathcal P}$ then ${H}\in {\mathcal F},$
thus $H\not\in {\rm Forb}_{\rm }({\mathcal F}).$
If ${ H}$  satisfies ${\mathcal P}$ but $H\not\in {\mathcal F}$ i.e. 
$H$ contains 
an $F\in {\mathcal F}$ as a subgraph, 
then, since ${\mathcal P}$ is hereditary, 
$F$ satisfies ${\mathcal P}$, but it contradicts $F\in {\mathcal F}$ 
with the definition of ${\mathcal F}$.
)
\par
We consider $\ell+2$ 
colors, called {\em black}, {\em white}$_i$($i\in [\ell]$) and {\em invisible}.
Set ${\rm BW}:=\{black, white_1, \cdots, white_\ell\}$ 
and ${\rm BI}:=\{black, invisible\}$.
A {\bf black-induced subgraph} of a $(k,{\rm BW})$-graph $H$ is 
a $(k,{\rm BI})$-graph obtained from $H$ by deleting some (if necessary) 
vertices 
and recoloring {\em all} the white edges and some (if necessary) 
black edges in the invisible color, 
where an edge is white if-and-only-if the color of the edge is white$_i$ for some $i$. 
A $(k, {\rm BW})$-property ${\mathcal P}$ is 
{\bf monotone} if-and-only-if there exists a $(k, {\rm BI})$-property 
${\mathcal F}$ such that ${\mathcal P}$ consists of all 
$(k, {\rm BW})$-graphs containing no copy of $F$ 
as a black-induced subgraph for 
any $F\in {\mathcal F}.$ 
We denote ${\mathcal P}={\rm Forb}_{\rm bi}({\mathcal F})$ where 
bi stands for black-induced.
Any monotone property is hereditary.
(This is easy to see. 
Define $(k, {\rm BW})$-property $\ang{\mathcal F}$ from ${\mathcal F}$ 
by putting to $\ang{\mathcal F}$ all the BW-colored hypergraphs $F'$
obtained from an $F\in {\mathcal F}$ by recoloring each invisible edge 
of $F$ in non-invisible colors (any way). Then ${\rm Forb_{bi}}({\mathcal F})
={\rm Forb_{}}(\ang{\mathcal F}).$)
\par
Here it is not hard to see that 
\begin{eqnarray*}
{\rm ex}(n, {\rm Forb_{bi}}({\mathcal F}))={\log_2{\ell+1\over \ell}
\over {n\choose k}}
\max\{ \mbox{\rm the number of black edges in } H |
H\in {\rm Forb_{bi}}(n, {\mathcal F})
\}+\log_2\ell
\end{eqnarray*}
where ${\rm Forb_{bi}}(n,{\mathcal F})$ denotes 
the family of $(k,{\rm BW})$-graphs on vertex set $[n]$ in 
${\rm Forb}_{bi}(n, {\mathcal F}).$
 (Hint: This basically follows from the correspondence between an 
$H\in {\rm Forb}_{\rm bi}(n, {\mathcal F})$ and 
a ${\mathcal P}$-good $(k,2^{\rm BW}\setminus\{\emptyset\})$-graph $H'$ 
on $[n]$ 
where, for any $e\in {[n]\choose k}$, (i) $H(e)=black$ iff 
$H'(e)={\rm BW}$ and (ii) $H(e)=white_i$ for some $i$ iff 
$H'(e)=\{white_1,\cdots,white_\ell\}.$
Note that $H\in {\rm Forb}_{\rm bi}(n, {\mathcal F})$ iff 
the corresponding $H'$ is 
a ${\mathcal P}$-good $(k,2^{\rm BW}\setminus\{\emptyset\})$-graph.
Further note that 
$\sum_e\log |H'(e)|=
(\# \mbox{\rm black edges in } H)\log (\ell+1) 
+(\# \mbox{\rm white edges in } H)
\log \ell 
=(\# \mbox{\rm black edges in } H)\log {\ell+1\over \ell}
+ {n\choose k}
\log \ell.
$
)
\par
Thus, if $k=2$ and $\ell=1$ (i.e. the case of ordinary graphs) then 
the famous Erd\H{o}s-Stone theorem \cite{ErSt} implies that 
\begin{eqnarray}
{\rm ex}(n, {\rm Forb_{bi}({\mathcal F})})
=
\min_{F\in {\mathcal F}}\brkt{1-{1\over \chi(F)-1}+o(1)},\label{uqq08}
\end{eqnarray}
where $\chi(F)$ is the chromatic number of $(2,{\rm BI})$-graph 
$F$. 
\subsection{A brief history of this research area}
As far as I know, 
all previous researchers have dealt with the case of two colors, 
black and white. We reset ${\rm BW}:=\{black, white\}.$
\subsubsection{Monotone properties for graphs}
Let $k=2.$
Erd\H{o}s, Kleitman and Rothschild \cite{EKR76} 
showed the theorem for $k=2$ and for ${\mathcal P}={\rm Forb_{bi}}(\{
K_\ell^{(2)}
\})$ where $K_\ell^{(2)}$ means the $\ell$-vertex $(2,\{black\})$-graphs 
with all edges black.
Using Szemer\'edi's regularity lemma, 
Erd\H{o}s, Frankl and R\"{o}dl \cite{EFR86} showed it 
for any monotone $(2,{\rm BW})$-property 
${\mathcal P}.$ 
\subsubsection{Hereditary properties for graphs}
Let $k=2$.
Pr\"{o}mel and Steger \cite{PS91,PS92,PS92b,PS93} began to study 
the hereditary property for ordinary graphs.
Pr\"{o}mel and Steger \cite{PS92b} showed our main theorem 
 for ${\mathcal P}={\rm Forb_{}}(\{F\})$ where 
$F$ is any fixed $(2,{\rm BW})$-graph.
Their proof has already used 
an early version of hypergraph regularity lemma, 
which was shown independently by 
Chung \cite{Chung91} and Steger \cite{Steger90}. 
(Their version partitions the vertex set only, without 
partitioning size-$2$ edges.)
Scheinermann and Zito \cite{SZ94} asked whether 
$\lim_{n\to \infty}{\log |{\rm Forb_{}}(n, {\mathcal F})|
\over {n\choose k}}$ exists where ${\mathcal F}$ is an arbitrary
 family of $(2,{\rm BW})$-graphs.
Answering this affirmatively, 
Alekseev \cite{A93} and Bollob\'as and Thomason \cite{BT97} 
independently showed 
our main theorem for $k=2$ and $|C|=2$. 
In this ordinary graph case,
 the definition of ${\rm ex}(n,{\mathcal P})$ can be restricted more.
Although we cannot use the Erd\H{o}s-Stone theorem for this case, 
these researchers showed that the limit takes a value from 
$0,{1\over 2}, {2\over 3}, {3\over 4},\cdots$ as in the monotone case 
(\ref{uqq08}).
\subsubsection{Monotone properties for hypergraphs}
Nagle and R\"{o}dl \cite{NR01} showed the theorem for $k=3$ 
and for ${\mathcal P}=
{\rm Forb_{bi}}(\{F\})$ where $F$ is a fixed $(3,{\rm BW})$-graph.
 Their proof method is based on 
 Frankl-R\"{o}dl's version of $(3,{\rm BW})$-graph regularity 
lemma (\cite{FR02}). 
Nagle, R\"{o}dl and Schacht \cite{NRS06} showed the theorem for general $k$ 
and for any monotone $(k,{\rm BW})$-property.
Their proof relies on their version of hypergraph regularity lemma 
\cite{RS,RS07r,RS07c}.
\subsubsection{Hereditary properties for hypergraphs}
Bollob\'as and Thomason \cite{BT95} 
showed the existence of $\lim_{n\to\infty} 
{\log |{\mathcal P}_n|\over {n\choose k}}
$  for general $k$ and $|C|=2$, without showing our main theorem.
Based on an extended Loomis-Whitney inequality \cite{LW49}, 
they showed that ${\log |{\mathcal P}_n|\over {n\choose k}}$ is non-increasing, 
which implies the existence of the limit.
(On the other hand, as mentioned previously, 
it is easy to see that ${\rm ex}(n,{\mathcal P})$ is non-increasing, and 
it will be also seen in our proof.)
When $k=3$ and $|C|=2,$ 
Kohayakawa et al. \cite{KNR} showed our main theorem, 
based on Frankl-R\"{o}dl's version of $(3,{\rm BW})$-graph regularity 
lemma.\par
We will consider the multicolor case instead of two-color case. 
Although this generalization has not yet been studied before,
it is iteself interesting and, 
furthermore, the colored hypergraphs are natural objects for 
regularity lemmas.
\subsection{Hypergraph regularity lemma}
The celebrated $(2,2)$-graph regularity lemma was discovered by Szemer\'edi 
\cite{Sz} as a lemma for 
his famous theorem on arithmetic progressions \cite{Sz75}.
Inspired by the success of the lemma in graph theory and others (see \cite{KSSS02}), 
research on quasi-random hypergraphs was initiated by 
Chung \cite{Chung,Chung91}, 
Chung-Graham \cite{CG,CG2,CG3}, 
Haviland-Thomason \cite{HT,HT2}, Steger \cite{Steger90} and 
Frankl-R\"{o}dl \cite{FR92}. 
For other early work, see \cite{BNS,CT}.
However these regularity lemmas are too weak for deep applications 
like the celebrated Szemer\'edi's progression theorem.
Frankl-R\"{o}dl \cite{FR02} suggested that if there 
exists a certain strong regularity lemma for $(k,2)$-graphs 
then it implies Szemer\'edi's 
theorem. 
They gave such a 
regularity lemma for $(3,2)$-graphs which implies Roth's theorem 
(i.e. the length-three case of Szemer\'edi's theorem).
(Also see \cite{NR03}.) 
Solymosi \cite{So03,So04} gave a short argument by which such a regularity lemma 
implies not only Szemer\'edi's theorem but also its multidimensional extension 
by Furstenberg-Katznelson \cite{FK78}.
\par
In 2003--2004, R\"{o}dl and his collaborators \cite{RSk04,NRS}
and Gowers \cite{G} independently obtained their $(k,2)$-graph
regularity lemmas which answers \cite{FR02}.
Slightly later, Tao \cite{Tao06} gave another
version. However while years have passed since their
preprints became available in the internet, applications of their
methods have been appearing more slowly than expected (R\"{o}dl
et al. \cite{RNSSK}). A major reason was that their methods are
rather cumbersome and technical for easy use in deeper
applications.
\par
It had been noted that unlike the situation for $(2,2)$-graphs, there are
several ways one might define regularity (i.e. a basic quasi-random 
property)for $(k,2)$-graphs
(R\"{o}dl-Skokan \cite[pp.1]{RSk04},Tao-Vu \cite[pp.455]{TV}).
Kohayakawa et al. \cite[pp.188]{KNR} say that
the basic objects involved in the Regularity Lemma and the
Counting Lemma are already somewhat technical and that simplifying
these lemmas would be of great interest. \par
The major purpose of this paper is to illustrate 
that a new regularity lemma \cite{I06} may meet 
these requirements. With \cite{I06}, we 
can naturally obtain strong quasi-random
properties not from one basic quasi-random property but from 
a simple construction of a certain partition. 
(Thus, the previous regularity lemmas correspond to 
our definition of partition, and the previous 
 counting lemmas correspond to our regularity lemma in our language.)
It gives a shorter elementary 
proof of Szemer\'edi's theorem as well as its multi-dimensional extension 
\cite{FK78}, with explicit density bounds.
It is achieved by a quite simple non-iterative (probabilistic)
construction which makes it easy to understand why it works. 
The construction of regularization is new even if we assume
we are working with $(2,2)$-graphs. Furthermore, it is strong; for
example, it generalizes edge-induced subgraph counting to multicolored 
vertex-induced subgraph counting, in the original setting itself.
\par
We have already seen two applications of \cite{I06} besides 
Szemer\'edi's theorem.
One of of the two is a positive 
answer \cite{I06m} to a question by Alon and Shapira 
 \cite{AS05} on property testing.
Even after R\"{o}dl et al. discovered their 
$(k,2)$-graph regularity lemma, 
they \cite{ARS} employed 
Frankl-R\"{o}dl's $(3,2)$-graph regularity lemma, 
instead of using their regularity lemma, 
and answered it for $k=3$.
Then wthiout developing the constructive argument due to \cite{AS05},
 they \cite{RSgene} 
answered it for general $k$ nonconstructively, relying on 
a non-constructive method of graph limits due to \cite{LS04,LS05}.
(See \cite{ES07} for hypergraph limits, which implies Frankl-R\"{o}dl\cite{FR92}'s 
preliminary regularity lemma.)
Independently from \cite{RSgene}, 
the constructive argument of \cite{AS05} was 
naturally extended by \cite{I06m} to general $k$ in the platform of 
 \cite{I06}. 
\par
The other example 
is a linear Ramsey number for bounded-degree hypergraphs.
Again 
after R\"{o}dl et al. discovered their regularity lemma,
Cooley et al. \cite{CFKO} and Nagle et al. \cite{NORS} independently 
obtained a linear Ramsey number for $k=3$, based on 
Frankl-R\"{o}dl's $(3,2)$-graph regularity lemma. 
Then Cooley et al. \cite{CFKO2} and Ishigami \cite{I06lr} independently 
obtained a linear Ramsey number for general $k.$ 
While both use the argument of \cite{CFKO},
they are based on R\"{o}dl-Schacht \cite{RS,RS07r,RS07c}'s 
$(k,2)$-graph regularity lemma (which is a variant of 
\cite{RSk04,NRS})
and on \cite{I06}'s $(k,C)$-graph regularity lemma, 
respectively. Repeating the argument of \cite{CFKO} 
is less cumbersome in the environment of \cite{I06}.
Furthermore, \cite{I06lr} deal with the multicolor case, while 
\cite{CFKO2} considers the two-color case only.
(Very recently, the mult-color result itself was reproved 
by Conlon et al.'s nice extension \cite{CFS} 
of Kostochka-R\"{o}dl's argument \cite{KR} 
with a significantly better bound and without any regularity lemmas.
But the techniques and the lemma in \cite{I06lr} are 
still worthwhile, and 
it would have some possibilities for some directions of its extensions, when hoping no 
good bounds.
In fact, we could have said the same about \cite{CRST}.)
\par
Here we will see the third example in this paper. 
We will see how easy we can 
extend the result by Kohayakawa et al. \cite{KNR} 
from $k=3$ to general $k.$
In fact, this is easier than the previous two examples.
However this may be a simple example quickly illustrating 
the way to apply the regularity lemma \cite{I06} and its potential, 
at least for readers who are not used to \cite{I06}.
Although 
the result of this paper
itself may be essentially obtained also by developing 
\cite{NRS06}, it would be cumbersome at least in the sense of 
Kohayakawa et al. \cite[pp.188]{KNR}, even for 
two-colored hypergraphs.
For multi-colored hypergraphs discussed here, 
it would be more cumbersome with thier environment.
\section{Statements of Regularity Lemma and Main Lemma}
In this paper, we denote by $\Prob$ and $\Ex$ the probability and expectation, 
respectively. We denote the conditional probability and exepctation by
$\Prob[\cdots|\cdots]$ and $\Ex[\cdots|\cdots].$ 
\begin{sett}\label{r0512}
Throughout this section, we fix a positive integer $r$ and 
an \lq index\rq\ set $\mathfrak{r}$ with $|\mathfrak{r}|=r.$ 
Also we fix a probability space 
$({\bf \Omega}_i,{\mathcal B}_i,\Prob)$ 
for each $i\in \mathfrak{r}$.
Assume that ${\bf \Omega}_i$ is finite (but its cardinality may not be 
constant) 
and ${\mathcal B}_i=2^{{\bf \Omega}_i}$ 
for the sake of simplicity.
Write ${\bf \Omega}:=({\bf \Omega}_i)_{i\in \mathfrak{r}}$.
\end{sett}
In order to avoid using technical words like mesurability or 
Fubini's theorem frequently to readers who are interested only in applications to 
discrete mathematics, 
we assume ${\bf \Omega}_i$ as a (non-empty) finite set.
However our argument should be extendable to a more general probability space.
For applications, ${\bf \Omega}_i$ would contain a huge number of vertices.
\par
For an integer $a$, we write $[a]:=\{1,2,\cdots,a\},$ and 
${\mathfrak{r}\choose [a]}:=\dot{\bigcup}_{i\in [a]}{\mathfrak{r}\choose i}
=\dot{\bigcup}_{i\in [a]}\{I\subset \mathfrak{r}| |I|=i\}.$
When $r$ sets $X_i, i\in {\mathfrak r},$ with indices from ${\mathfrak r}$ are 
called {\bf vertex sets}, 
we write $X_J:=\{Y\subset \dot{\bigcup}_{i\in J}X_i| |Y\cap X_j|=1 \forall j\in J\}$ 
whenever $J\subset {\mathfrak r}$.
\begin{df}[(Bound colored hyper)graphs]
Suppose Setup \ref{r0512}. 
A {\bf $k$-bound $(b_i)_{i\in [k]}$-colored ($\mathfrak{r}$-partite hyper)graph} $H$ 
is a triple $((X_i)_{i\in\mathfrak{r}},({C}_I)_{I\in {\mathfrak{r}\choose [k]}},
(\gamma_I)_{I\in {\mathfrak{r}\choose [k]}}
)$ where (1) each $X_i$ is a set called a \lq vertex set,\rq\ (2)
${C}_I$ is a set with at most $b_{|I|}$ elements, and 
(3) $\gamma_I$ is a map from $X_I$ to ${C}_I.$
We write $V(H)=\dot{\bigcup}_{i\in \mathfrak{r}}
V_i(H)=\dot{\bigcup}_{i\in \mathfrak{r}}X_i
$ and 
${\rm C}_I(H)={C}_I.$ 
Each element of $V(H)$ is called a {\bf vertex}.
Each element $e\in V_I(H)=X_I, I\in {\mathfrak{r}\choose [k]},$ is called 
 an {\bf (index-$I$ size-$|I|$) edge}.
Each member in ${\rm C}_I(H)$ is a {\bf (face-)color (of index $I$)}.
Write $H(e)=\gamma_I(e)$ for each $I.$ 
Put ${\rm C}_i(H):=\dot{\bigcup}_{I\in {\mathfrak{r}\choose i}}{\rm C}_I(H)$ 
for $i\in [k].$
\par
Let $I\in {\mathfrak{r}\choose [k]}$ 
and $e\in V_I(H).$ 
For another index $\emptyset\not=
J\subset I$, we denote by $e|_J$ the index-$J$ edge 
$e\setminus \brkt{\bigcup_{j\in I\setminus J}X_j}\in V_J(H)$.
We define the {\bf frame-color}
and {\bf total-color} of $e$ by ${H}(\partial
e):=({H}(e|_J)|\,\emptyset\not=J\subsetneq I)$ and by
${H}(\ang{e})=H\ang{e}:=({H}(e|_J)|\,\emptyset\not=J\subsetneq I).$
Write ${\rm TC}_I(H):=\{H\ang{e}|\,{e}\in X_I\},$ ${\rm
TC}_s(H):=\bigcup_{I\in {\mathfrak{r}\choose s}}{\rm TC}_I(H),$
and ${\rm TC}(H):=\bigcup_{s\in  [k]}{\rm TC}_s(H).$
\par
A {\bf ($k$-bound) (simplicial-)complex} is a $k$-bound
(colored ${\mathfrak r}$-partite hyper)graph 
such that for each $I\in {\mathfrak{r}\choose [k]}$ 
there exists at most one index-$I$ color called \lq invisible\rq\ 
and that if (the color of) an edge $e$ is invisible then 
any edge $e^*\supset e$ is invisible. An edge or its color 
is {\bf visible} if it is not invisible. 
\par
For a $k$-bound graph ${\bf G}$ on ${\bf \Omega}$ and $s\le k$, 
let ${\mathcal S}_{r,s,h,{\bf G}}={\mathcal S}_{s,h,{\bf G}}$ be the set of $s$-bound 
simplicial-complexes $S$ such that 
(1) each of the $r$ vertex sets contains exactly $h$ vertices and 
that (2)
for any $I\in {\mathfrak{r}\choose [s]}$ 
 there is an injection from the index-$I$ visible colors of $S$ to the 
index-$I$ colors of ${\bf G}$.
(When a visible color $\mathfrak{c}$ of $S$ corresponds to another color $\mathfrak{c}'$ of 
${\bf G}$, we simply write $\mathfrak{c}=\mathfrak{c}'$ 
without presenting the injection explicitly.)
For $S\in {\mathcal S}_{s,h,{\bf G}}$, we denote by ${\mathbb V}_I(S)$ the set of 
index-$I$ visible edges. Write ${\mathbb V}_i(S):=\bigcup_{I\in {\mathfrak{r}\choose 
i}}{\mathbb V}_I(S)$ and ${\mathbb V}(S):=\bigcup_i {\mathbb V}_i(S).$
\end{df}
\begin{df}[Partitionwise maps]
A {\bf partitionwise map} $\varphi$ is 
a map 
 from $r$ vertex sets $W_i,i\in \mathfrak{r},$ with $|W_i|<\infty$ to {\em the} 
$r$ vertex sets (probability spaces)$
U_i,i\in\mathfrak{r}$,
such that 
each $w\in W_i$ is mappped into $U_i.$
We denote by $\Phi((W_i)_{i\in\mathfrak{r}}, (U_i)_{i\in\mathfrak{r}})$ 
or $\Phi(\bigcup_{i\in\mathfrak{r}}W_i,\bigcup_{i\in\mathfrak{r}}U_i)$ 
the set of partitionwise maps from $(W_i)_i$ to $(U_i)_i$.
If $U_i={\bf \Omega}_i$ or $U_i$ is obvious then we omit them.
A partitionwise map is {\bf random} if-and-only-if 
each $w\in W_i$ is mutually-independently mapped at random 
according to the probability space 
${\bf \Omega}_i$.
\end{df}
We define the regularity of hypergraphs.
\begin{df}[Regularity]
Let ${\bf G}$ be a $k$-bound graph on $
{\bf \Omega}$.
For 
$\vec{\mathfrak{c}}=(\mathfrak{c}_J)_{J\subset I}\in {\rm TC}_I({\bf G}), 
I\in {\mathfrak{r}\choose [k]}$, we define {\bf relative density}
\begin{eqnarray*}
{\bf d}_{\bf G}(\vec{\mathfrak{c}}):=
\Prob_{{\bf e}\in {\bf \Omega}_I
}[
{\bf G}({\bf e})=\mathfrak{c}_I
|
{\bf G}(\partial{\bf e})=
(\mathfrak{c}_J)_{J\subsetneq I}
].
\end{eqnarray*}
\par
For a positive integer $h$ and a real 
$\epsilon>0$, we say that 
 ${\bf G}$ is 
 {\bf $(\epsilon,h)$-regular} 
if-and-only-if 
there exists a function
${ \delta}: {\rm TC}({\bf G})\to [0,\infty)$ 
such that 
\begin{eqnarray}
\hspace{-5mm}
{\rm (i)}&
 \Prob_{\phi\in\Phi(V(S))}
[
{\bf G}(\phi(e))=S(e)
\forall e\in\mathbb{V}(S)
]=
\displaystyle\prod_{e\in {\mathbb V}(S)}
\brkt{
{\bf d}_{\bf G}(S\ang{e})
\dot{\pm}
\delta(S\ang{e})
}&
\forall S\in {\mathcal S}_{k,h,{\bf G}},
\label{a0719}
\\
\hspace{-5mm}
{\rm (ii)}& \Ex_{{\bf e}\in {\bf \Omega}_I}[\delta({\bf G}\ang{\bf e})]\le 
\epsilon/
|{\rm C}_I({\bf G})|
& \forall I\in {\mathfrak{r}\choose [k]},
\label{a0726}
\end{eqnarray}
where $a\dot{\pm}b$ means (the interval of) numbers $c$ with 
$\max\{0,a-b\}\le c\le \min\{1,a+b\}$.
\par
A {\bf subdivision }  of a $k$-bound graph ${\bf G}$ on ${\bf \Omega}$ 
is a $k$-bound graph ${\bf G}^*$ on the same ${\bf \Omega}$ such that 
\\
(i) for any size-$k$ edge ${\bf e}\in {\bf \Omega}_I$ with 
$I\in {\mathfrak{r}\choose k},$ 
it holds that ${\bf G}^*({\bf e})={\bf G}({\bf e}),$ and \\
(ii) for any two edges ${\bf e}, {\bf e'}\in {\bf \Omega}_I$ with 
$I\in {\mathfrak{r}\choose [k-1]}$, 
if ${\bf G}^*({\bf e})={\bf G}^*({\bf e'})$ then ${\bf G}({\bf e})={\bf G}({\bf e'}).$
\end{df}
\begin{tha}[Hypergraph Regularity Lemma in \cite{I06}] 
\label{a0721}
Let $r\ge k,h,\vec{b}=(b_i)_{i\in [k]}$ be positive integers, and 
$\epsilon>0
$ a real.
Then there exist integers 
$\widetilde{b}_1\ge \cdots\ge \widetilde{b}_{k-1}
$ (independent from ${\bf \Omega}$)
such that if ${\bf G}$  is 
a $\vec{b}$-colored ($k$-bound $r$-partite hyper)graph on ${\bf \Omega}$ then 
there exists an $(\epsilon,h)$-regular $(\widetilde{b}_1,\cdots,
\widetilde{b}_{k-1},b_k)$-colored 
subdivision ${\bf G}^*$ of ${\bf G}.$ 
\end{tha}
For the simple way to construct such a subdivision in 
Theorem \ref{a0721}, see  \cite{I06}.
\par
I think that 
readers who understand our version of regularity lemma 
will feel that the proof in the next section 
is easy, once the statement is given.
\section{Proof of the Theorem}
$\bullet$ 
While a $(k,C)$-graph $H$ can be expressed as a function 
$H: {V(H)\choose k}\to C$, 
an {\bf $r$-partite} $(k,C)$-graph is a function 
expressed as $H: \{e\in {V(H)\choose k}
 ; |e\cap W_i|\le 1 \forall i\in [r]
\}
\to C$ 
for some vertex partition $V(H)=W_1\dot{\cup}\cdots\dot{\cup}W_r.$ 
That is, every edge considered in $r$-partite $(k,C)$-graph $H$ is 
\lq partitionwise\rq\, i.e. 
every edge contains at most one vertex in each $W_i.$
\par
Let $b_k$ be a positive integer.
Let $C$ be a color class with $|C|=b_k$ and ${\mathcal P}$ be 
a $(k,C)$-property.
\par
It is easy to see that 
\begin{eqnarray}
\log_2 
|{\mathcal P}_n|\ge 
{{n\choose k}{\rm ex}(n, {\mathcal P})}.\label{uqq03a}
\end{eqnarray}
In fact, suppose an $n$-vertex ${\mathcal P}$-good 
$(k,2^C\setminus\{\emptyset\})$-graph ${H}$ with 
${n\choose k}^{-1}\sum_{e\in {V(H)\choose k}}\log |{H}({e})|=
{\rm ex}(n, {\mathcal P})$.
No matter how each edge $e\in {V(H)\choose k}$ is recolored
 by a color in $H(e)$, 
the resulting $(k,C)$-graph satisfies ${\mathcal P}$. 
The number of such $(k,C)$-graphs is 
$\prod_e |H(e)|=2^{\sum_e\log_2|H(e)|}
=2^{{n\choose k}
{\rm ex}(n, {\mathcal P})},$ 
yielding  (\ref{uqq03a}).
\par
Thus our goal is to show for any constant $\eta>0,$ 
\begin{eqnarray}
\log_2 |{\mathcal P}_n|\le {{n\choose k}
({\rm ex}(n, {\mathcal P})+1.1\eta)} .\label{uqq03b}
\end{eqnarray}
We set the following parameters
\begin{eqnarray}
k, b_k, 1/\eta \ll r\ll 1/\alpha \ll n\label{uqq06}
\end{eqnarray}
where $r,\alpha,n$ depend also on ${\mathcal P}$.
\par
Let ${\bf G}\in {\mathcal P}_n.$ 
We set the vertex set $[n]={\bf \Omega}_1\dot{\cup}\cdots\dot{\cup}{\bf \Omega}_r$ 
with $N:=|{\bf \Omega}_i|=n/r$ so that 
${\bf \Omega}_i:=\{(i-1)N+1,\cdots,iN\}$ for each $i\in [r]$.
Here we assume that $r$ divides $n.$
If not, we remove at most $r-1$ vertices so that all partite sets have 
 the same vertices. 
Note that the resulting $(k,C)$-graph 
still satisfies the property ${\mathcal P}$ 
due to the definition of hereditary.
The number of possible color patterns 
of edges containing the removed vertices is at most 
\begin{eqnarray}
b_k^{(r-1){n-1\choose k-1}}=2^{O(n^{k-1})}.\label{uqq04}
\end{eqnarray}
From now on, we will never look at any \lq non-partitionwise\rq\ edge.
That is, in the resulting $r$-partite $(k,C)$-graph, 
any edge contains at most one vertex in each partite set ${\bf \Omega}_i$.
The number of the non-partitionwise edges 
is at most $(k-1)r^{k-1}N^k=(k-1)r^{k-1} (n/r)^k=(k-1)n^k/r$. 
Hence the number of possible color patterns of the non-partitionwise edges is at most 
\begin{eqnarray}
b_k^{kn^k/r}.
\label{uqq04b}
\end{eqnarray}
\par
Next we color in \lq white\rq\ 
all the edges ${\bf e}\in {\bf \Omega}_I$ of size at most $k-1$, 
i.e. $I\in {[r]\choose [k-1]}$.
For this resulting $r$-partite $k$-bound $(1,\cdots,1,b_k)$-colored 
graph, we apply 
 the regularity lemma (Theorem \ref{a0721}) 
with $r:=r,k:=k,h:=1,\vec{b}:=(1,\cdots,1,b_k)$ and with
\begin{eqnarray}
\epsilon:=
\brkt{{\alpha\over 11\cdot 2^k}}^2
\label{xx17a0}
\end{eqnarray}
 and obtain an 
$(\epsilon,1)$-regular subdivision 
${\bf G}^*$ which is $(\widetilde{b}_1,\cdots,
\widetilde{b}_{k-1},\widetilde{b}_k=b_k)$-colored where 
\begin{eqnarray*}
k,r,b_k,1/\alpha \ll 
\widetilde{b}_{k-1},\cdots,\widetilde{b}_{1}\ll n.
\label{xx17a}
\end{eqnarray*}
The number of possible color patterns of edges of size at most $k-1$ 
is at most 
\begin{eqnarray}
\widetilde{b}_1^{rN}\widetilde{b}_2^{{r\choose 2}N^2}
\cdots \widetilde{b}_{k-1}^{{r\choose k-1}N^{k-1}}
=2^{O(n^{k-1})}.
\label{uqq04d}
\end{eqnarray}
A total color $\vec{\mathfrak{c}}=(\mathfrak{c}_J)_{J\subset I}
\in {\rm TC}_I({\bf G^*})$ 
 with 
$I\in {[r]\choose [k]}$ is called {\bf exceptional} 
if-and-only-if there exists $I'\subset I$ such that 
${\bf d}_{{\bf G}^*}((\mathfrak{c}_J)_{J\subset I'})<
\sqrt{\epsilon}
/|{\rm C}_{I'}({\bf G^*})|
$ or 
$\delta((\mathfrak{c}_J)_{J\subset I'})>
0.1\sqrt{\epsilon}/|{\rm C}_{I'}({\bf G^*})|
$ where 
$\delta(\cdot)$ is a function associated with ${\bf G}^*$.
An edge is said to be exceptional if-and-only-if
 its total color 
${\bf G^*}\ang{\bf e}$ is 
exceptional.
For any index $I$, it easily follows that 
\begin{eqnarray}
&&
\Prob_{{\bf e}\in {\bf \Omega}_I}[{\bf G}^*\ang{\bf e}
\mbox{ is exceptional}]\nonumber
\\
&\le &
\sum_{J\subset I}
\brkt{
\Prob_{{\bf e}\in {\bf \Omega}_{J}}[
{\bf d}_{{\bf G^*}}({\bf G^*}\ang{\bf e})
<{\sqrt{\epsilon}\over |{\rm C}_{J}({\bf G^*})|}]
+
\Prob_{{\bf e}\in {\bf \Omega}_{J}}
[\delta({\bf G^*}\ang{\bf e})
>{0.1\sqrt{\epsilon}\over |{\rm C}_{J}({\bf G^*})|}]
}\nonumber
\\
&\stackrel{(**)}{\le}_{(\ref{a0726})}
 &
\sum_{J\subset I}
\brkt{
\sqrt{\epsilon}
+
{
\epsilon\over
0.1\sqrt{\epsilon}
}
}
<
11\cdot 2^k\sqrt{\epsilon}
\stackrel{(\ref{xx17a0})}{=}
\alpha\label{uqq06e}
\end{eqnarray}
where in the above (**) we used the fact that 
\begin{eqnarray*}
&&
\Prob_{{\bf e}\in {\bf \Omega}_J}\left[
\Prob_{{\bf e'}\in {\bf \Omega}_J}
[{\bf G^*}({\bf e'})
={\bf G^*}({\bf e})
|{\bf e'}\stackrel{\partial {\bf G^*}}{\approx}{\bf e}
]
\le
{\sqrt{\epsilon}\over |{\rm C}_J({\bf G^*})|}
\right]
\nonumber
\\
&=&\sum_{\mathfrak{c}_J\in {\rm C}_J({\bf G^*})}
\Prob_{{\bf e}\in {\bf \Omega}_J}\left[
{\bf G^*}({\bf e})=\mathfrak{c}_J \mbox{ and }
\Prob_{{\bf e'}\in {\bf \Omega}_J}
[{\bf G^*}({\bf e'})
=\mathfrak{c}_J
|{\bf e'}\stackrel{\partial {\bf G^*}}{\approx}{\bf e}
]
\le
{\sqrt{\epsilon}\over |{\rm C}_J({\bf G^*})|}
\right]\nonumber
\\
&\le &\sum_{\mathfrak{c}_J\in {\rm C}_J({\bf G^*})}
1\cdot
\Prob_{{\bf e}\in {\bf \Omega}_J}\left[
{\bf G^*}({\bf e})=\mathfrak{c}_J \left|
\Prob_{{\bf e'}\in {\bf \Omega}_J}
[{\bf G^*}({\bf e'})
=\mathfrak{c}_J
|{\bf e'}\stackrel{\partial {\bf G^*}}{\approx}{\bf e}
]
\le
{\sqrt{\epsilon}\over |{\rm C}_J({\bf G^*})|}
\right.
\right]\nonumber
\\
&= &\sum_{\mathfrak{c}_J\in {\rm C}_J({\bf G^*})}
\Ex_{{\bf e}\in {\bf \Omega}_J}\left[
\Prob_{{\bf \tilde{e}}\in {\bf \Omega}_J}
[
{\bf G^*}({\bf \tilde{e}})=\mathfrak{c}_J 
|
{\bf \tilde{e}}\stackrel{\partial {\bf G^*}}{\approx}
{\bf e}
]
\left|
\Prob_{{\bf e'}\in {\bf \Omega}_J}
[{\bf G^*}({\bf e'})
=\mathfrak{c}_J
|{\bf e'}\stackrel{\partial {\bf G^*}}{\approx}{\bf e}
]
\le
{\sqrt{\epsilon}\over |{\rm C}_J({\bf G^*})|}
\right.
\right]\nonumber
\\
&& \hspace{1.5in}\mbox{($\because$ the conditional part depends only on ${\bf G^*}(\partial {\bf e})$)}
\nonumber
\\
&\le  &\sum_{\mathfrak{c}_J\in {\rm C}_J({\bf G^*})}
\Ex_{{\bf e}\in {\bf \Omega}_J}\left[
{\sqrt{\epsilon}\over |{\rm C}_J({\bf G^*})|}
\left|
\Prob_{{\bf e'}\in {\bf \Omega}_J}
[{\bf G^*}({\bf e'})
=\mathfrak{c}_J
|{\bf e'}\stackrel{\partial {\bf G^*}}{\approx}{\bf e}
]
\le
{\sqrt{\epsilon}\over |{\rm C}_J({\bf G^*})|}
\right.
\right]
=\sqrt{\epsilon}.
\end{eqnarray*}
Thus, since 
Stirling inequality implies ${a\choose b}
\le a^b/\sqrt{2\pi b}(b/{\rm e})^b
< (a{\rm e}/b)^b$ with ${\rm e}:=2.71828\cdots$,
 the number of possible distributions of exceptional size-$k$ edges 
in one $I\in {[r]\choose k}$ is at most 
\begin{eqnarray*}
\sum_{j\le \alpha N^k}
{N^k\choose j}\le \alpha N^k
({\rm e}/\alpha)^{\alpha N^k}=
2^{\alpha\log_2 ({\rm e}/\alpha)\cdot N^k+o(N)}
=2^{\alpha\log_2 ({\rm e}/\alpha)\cdot (n/r)^k+o(n)}
\end{eqnarray*}
when ignoring the face-colors of the exceptional edges.
If we count the face-colors of the 
exceptional edges in all $I\in {[r]\choose k}$, 
it becomes at most 
\begin{eqnarray}
\brkt{
2^{\alpha\log_2 
({\rm e}/\alpha)\cdot (n/r)^k+o(n)}\cdot b_k^{\alpha (n/r)^k}
}^{r\choose k}
.
\label{uqq04e}
\end{eqnarray}
We are now defining an 
\lq almost ${\mathcal P}$-good\rq\ $r$-partite 
$(k,2^C\setminus\{\emptyset\})$-graph ${\bf H}$ as follows: 
The vertex set of ${\bf H}$ is the same as ${\bf G^*}.$ 
For any size-$k$ edge ${\bf e}\in {\bf \Omega}_I$ with $I\in {[r]\choose k},$ 
if ${\bf e}$ is not exceptional then 
we assign ${\bf e}$ the set of all 
face-colors $\mathfrak{c}_I\in {\rm C}_I({\bf G^*})
={\rm C}_I({\bf G})$ such that $(\mathfrak{c}_J)_{J\subset I}$ 
with $(\mathfrak{c}_J)_{J\subsetneq I}:={\bf G^*}(\partial {\bf e})$ is 
not exceptional. (Note that such a face-color $\mathfrak{c}_I$ exists.)
If ${\bf e}$ is exceptional, then we do not recolor ${\bf e}.$ 
We remove all edges of size at most $k-1$. The resulting $r$-partite 
$(k,2^C\setminus\{\emptyset\})$-graph is denoted by ${\bf H}.$
\begin{claim}\label{uqq04g}
$\sum_{I\in {[r]\choose k}}\sum_{{\bf e}\in {\bf \Omega}_I}
\log_2 |{\bf H}({\bf e})|\le 
({\rm ex}(n, {\mathcal P})+\eta) {n\choose k}.$
\end{claim}
Before proving this claim, we show that it implies our main theorem.
\par
By (\ref{uqq04e}), when 
 we fix the colors of the edges of size at most $k-1$ in ${\bf G^*},$ 
the number of possible face-color patterns of the size-$k$ edges 
in ${\bf H}$ is at most 
\begin{eqnarray}
\brkt{
2^{\alpha\log_2 
({\rm e}/\alpha)\cdot (n/r)^k+o(n)}\cdot b_k^{\alpha (n/r)^k}
}^{r\choose k}
\cdot 
(2^{b_k}-1)^{
{r\choose k}\widetilde{b}_{k-1}^{k\choose k-1}\cdots \widetilde{b}_1^{k\choose 1}
}\le 
2^{\alpha\log_2 ({\rm e}/\alpha)\cdot n^k/k!}\cdot b_k^{\alpha n^k/k!}
,\label{uqq06b}
\end{eqnarray}
where the exponent $
{r\choose k}\widetilde{b}_{k-1}^{k\choose k-1}\cdots \widetilde{b}_1^{k\choose 1}
$ means the number of possible frame-colors of size $k.$
\par
Once we fix the colors of the edges of size at most $k-1$ 
in ${\bf G^*}$ and fix 
the colors of the size-$k$ edges in ${\bf H},$ 
the number of possible graphs ${\bf G^*}$ producing 
${\bf H}$ is at most 
\begin{eqnarray}
\prod_{I\in {[r]\choose k}}
\prod_{{\bf e}\in {\bf \Omega}_I}
|{\bf H}({\bf e})|=
2^{\sum_{I,{\bf e}}\log_2 |{\bf H}({\bf e})|}
\stackrel{Claim \ref{uqq04g}}{\le}
 2^{{n\choose k}(\eta+
{\rm ex}(n,{\mathcal P}))}.\label{uqq06a}
\end{eqnarray}
Finally, by (\ref{uqq04}), (\ref{uqq04b}), 
(\ref{uqq04d}), 
(\ref{uqq06b}) and (\ref{uqq06a}), 
the number of graphs ${\bf G}$ in ${\mathcal P}_n$ is at most 
\begin{eqnarray*}
&&
2^{O(n^{k-1})}
b_k^{kn^k/r}2^{O(n^{k-1})}
2^{(\log_2({\rm e}/\alpha)+\log_2 b_k)\alpha n^k/k!}
2^{{n\choose k}(\eta+{\rm ex}(n,{\mathcal P}))}
\nonumber
\\
&=&
2^{
({\rm ex}(n,{\mathcal P})
+\eta){n\choose k}+({k\over r}\log_2 b_k)n^k
+(\log_2({\rm e}/\alpha)+\log_2 b_k)\alpha n^k/k!+O(n^{k-1})
}\nonumber
\\
&
\stackrel{(\ref{uqq06})}{\le}&
2^{
({\rm ex}(n, {\mathcal P})+1.1\eta) {n\choose k}},
\end{eqnarray*}
showing (\ref{uqq03b}). This together with (\ref{uqq03a}) 
completes the proof of Theorem \ref{uqq03}, if Claim \ref{uqq04g} is true.
\qed
\bigskip
\\
{\bf Proof of Claim \ref{uqq04g} :}
Assume that 
\begin{eqnarray}
\sum_{I\in {[r]\choose k}}\sum_{{\bf e}\in {\bf \Omega}_I}\log |{\bf H}({\bf e})|
>({\rm ex}(n, {\mathcal P})+\eta){n\choose k}.
\label{uqq06f}
\end{eqnarray}
\par
It is easily seen that the value 
${\rm ex}(\ell, {\mathcal P})\ge 0$ 
is non-increasing when $\ell$ is growing.
(Indeed, suppose ${\rm ex}(\ell+1,{\mathcal P})
>
{\rm ex}(\ell,{\mathcal P}).$ 
Then there exists an $(\ell+1)$-vertex ${\mathcal P}$-good 
 $(k,2^C\setminus\{\emptyset\})$-graph $F$ on 
vertex set $[\ell+1]$
with $
\Ex_{e\in {[\ell+1]\choose k}}\log |F(e)|
={\rm ex}(\ell+1,{\mathcal P}).
$ 
Since $\Ex_u\Ex_{e\in {[\ell+1]\setminus\{u\}\choose k}}
=\Ex_{e\in {[\ell+1]\choose k}},
$
there exists a vertex $u\in [\ell+1]$ such that, by deleting $u$ from $F,$
 the resulting ${\mathcal P}$-good graph $F^-$ satisfies the property that 
$\Ex_{e\in {[\ell+1]\setminus\{u\}\choose k }} \log|F^-(e)|\ge 
\Ex_{e\in {[\ell+1]\choose k}}\log |F(e)|
=
{\rm ex}(\ell+1, {\mathcal P})
>
{\rm ex}(\ell, {\mathcal P})$ by our assumption, 
contradicting definition of ${\rm ex}(\ell, {\mathcal P}).$
)
\par
Since the above guarantees the existence of $\lim_{\ell\to\infty}
{\rm ex}(\ell,{\mathcal P})\in [0,1], 
$
we can take an $r$ so that 
\begin{eqnarray}
{\rm ex}(r,{\mathcal P})
\ge 
{\rm ex}(n,{\mathcal P})
\ge 
{\rm ex}(r,{\mathcal P})
-0.1\eta\label{uqq05}
\end{eqnarray}
where 
$
k,1/\eta, {\mathcal P}\ll 
 r\le n
$ by (\ref{uqq06}).
\par 
Take 
${\bf e_0}\in {\bf \Omega}_{[r]}$ randomly.
We have 
\begin{eqnarray}
\Prob_{{\bf e_0}\in {\bf \Omega}_{[r]}}\left[
\exists {\bf e}\in {{\bf e_0} \choose k}, 
{\bf G^*}\ang{\bf e} \mbox{ is exceptional}
\right]
\stackrel{(\ref{uqq06e})}{<}{r\choose k}\alpha<1
\label{uqq06h}
\end{eqnarray}
where 
$
k,r\ll 1/\alpha
$ by (\ref{uqq06}).
On the other hand, we have that 
\begin{eqnarray}
\Ex_{{\bf e_0}\in {\bf \Omega}_{[r]}}
[\Ex_{{\bf e}\in {{\bf e_0}\choose k}}
\log |{\bf H}({\bf e})|
]
&=&
\Ex_{I\in {[r]\choose k}}
\Ex_{{\bf e}\in {\bf \Omega}_I}[\log |{\bf H}({\bf e})|
]
\nonumber\\
&=&{
1
\over {r\choose k}N^k
}
\sum_{I\in {[r]\choose k}}
\sum_{{\bf e}\in {\bf \Omega}_I}\log |{\bf H}({\bf e})|
\nonumber\\
&\stackrel{(\ref{uqq06f})}{\ge} &
{
{n\choose k}({\rm ex}(n, {\mathcal P})+\eta)
\over {r\choose k}(n/r)^k}
\nonumber\\
&\stackrel{(\ref{uqq05})}{\ge}
 &
{(1-o(1)){(1/k!)}(
{\rm ex}(r,{\mathcal P})+0.9\eta)
\over {r\choose k}(1/r)^k}
\nonumber\\
&> &
 {\rm ex}(r,{\mathcal P})
+0.9\eta.\label{uqq06g}
\end{eqnarray}
Therefore by (\ref{uqq06h}) and (\ref{uqq06g}) 
there exists an ${\bf e}_0\in {\bf \Omega}_{[r]}$
 such that (i) every ${\bf e}\in {{\bf e_0}\choose k}$ is not exceptional 
and that (ii) \begin{eqnarray}
\Ex_{{\bf e}\in {{\bf e_0}\choose k}}\log |{\bf H}({\bf e})|
&>& {{\rm ex}(r, 
{\mathcal P})+0.9\eta-{r\choose k}\alpha\cdot\log b_k
\over 1-{r\choose k}\alpha
}
\nonumber
\\
&>& ({\rm ex}(r, 
{\mathcal P})+0.9\eta-{r\choose k}\alpha\log b_k
)( 1+{r\choose k}\alpha)
\nonumber
\\
&>& 
{\rm ex}(r, {\mathcal P})+0.9\eta-{r\choose k}\alpha\log b_k
-{r\choose k}^2\alpha^2\log b_k\nonumber
\\
&>& 
{\rm ex}(r, {\mathcal P})+0.8\eta\label{uqq06i}
\end{eqnarray}
where we used $|{\bf H}({\bf e})|\le b_k$ in the first inequality 
and where we used in the last inequality 
$ k,b_k,1/\eta,r\ll 1/\alpha$ because of (\ref{uqq06}).
Due to the definition of ${\rm ex}(r, {\mathcal P})$ with (\ref{uqq06i}), 
we can make an $r$-vertex $(k,C)$-graph 
$H\not\in {\mathcal P}$ (on vertex set $
V(H):={\bf e_0}=\{{\bf v}_1,\cdots,
{\bf v}_r\}$) from ${\bf e_0}$
by 
assigning each ${\bf e}\in {{\bf e_0}\choose k}$
a color $\mathfrak{c}\in {\bf H}({\bf e})\subset {\rm C}_k({\bf G})$.
Furthermore we make a simplicial-complex $S\in {\mathcal S}_{k,h,{\bf G^*}}$ 
from this $H$ by assigning each ${\bf e}\in {{\bf e_0}\choose [k-1]}$ 
a color in ${\bf G^*}({\bf e}).$ Here all edges of $S$ are visible.
Since all edges of $S$ are not exceptional by (i), 
it follows from Theorem \ref{a0721} and from 
the definition of exceptional total-colors that 
\begin{eqnarray}
\Prob_{\phi\in \Phi(V(H))}
[{\bf G^*}(\phi(e))=H(e) \forall e\in {V(H)\choose k}
]\nonumber
&\ge&
\Prob_{\phi\in \Phi(V(S))}
[{\bf G^*}(\phi(e))=S(e) \forall e\in 
{V(H)\choose [k]}
]\nonumber
\\
&\stackrel{(\ref{a0719})}{\ge}&
 \prod_{e\in \mathbb{V}(S)}({\bf d}_{\bf G^*}(S\ang{e})-\delta(S\ang{e}))
\nonumber\\
&\ge& \prod_{e\in \mathbb{V}(S)}0.9
{\sqrt{\epsilon}\over \widetilde{b}_{|e|}}>0.\label{uqq04f}\nonumber
\end{eqnarray}
Therefore 
$H$ can be a subgraph of the $(k,C)$-graph ${\bf G}$ by relabeling 
the vertices, 
since the color of any size-$k$ edge was not recolored when regularizing ${\bf G}.$
Since ${\mathcal P}$ is hereditary and contains ${\bf G}$ by our assumption, 
any subgraph of ${\bf G}$  satisfies ${\mathcal P}$. Thus 
$H\in {\mathcal P}$. However 
the definition of $H$ implies that $H\not\in {\mathcal P}.$
This contradiction completes the proof of Claim \ref{uqq04g}.
\claimqed

\end{document}